\documentclass[a4paper,11pt,
twoside]{article}

\usepackage{amsmath, amsthm, amssymb,bm}

\numberwithin{equation}{section}

\newcommand{\bC}{\mathbb{C}}

\newcommand{\bQ}{\mathbb{Q}}
\newcommand{\bT}{\mathbb{T}}

\newcommand{\bZ}{\mathbb{Z}}

\newcommand{\mf}[1]{\mathfrak{#1}}

\def\pf{\mathrm{pf}}

\def\sgn{\mathrm{sgn}}

\newcommand{\hdet}[1]{{\textstyle \det^{[#1]} }}
\newcommand{\hper}[1]{{\textstyle \mathrm{per}^{[#1]} }}
\newcommand{\hpf}[1]{{\textstyle \mathrm{pf}^{[#1]} }}
\newcommand{\hPf}[1]{{\textstyle \mathrm{Pf}^{[#1]} }}

\def\dd{\mathrm{d}}

\def\bx{\mathbf{x}}
\def\by{\mathbf{y}}

\def\({ \left( }
\def\){ \right)}
\def\[{ \left[ }
\def\]{ \right]}

\theoremstyle{plain}
\newtheorem{thm}{Theorem}[section]
\newtheorem{prop}[thm]{Proposition}
\newtheorem{lem}[thm]{Lemma}
\newtheorem{cor}[thm]{Corollary}
\theoremstyle{definition}

\newtheorem{example}{Example}[section]
\newtheorem{remark}{Remark}[section]
\theoremstyle{conjecture}

\theoremstyle{problem}

\title{\bfseries Hyperdeterminantal expressions for Jack functions of rectangular shapes}
\author{\textsc{Sho MATSUMOTO}
\thanks{Research Fellow of the Japan Society 
for the Promotion of Science, partially supported by Grant-in-Aid 
for Scientific Research (C) No. 17006193.}%
\\
{\it \small Faculty of Mathematics, Kyushu University.} \\
{\it \small Hakozaki, Higashi-ku, Fukuoka, 812-8581 Japan.} \\
\texttt{shom@math.kyushu-u.ac.jp}
}

\date{\empty}

\begin{document}
\maketitle

\begin{abstract}
We derive a Jacobi-Trudi type formula for Jack functions of rectangular shapes.
In this formula, we make use of a hyperdeterminant, 
which is Cayley's simple generalization of the determinant. 
In addition,
after developing the general theory of hyperdeterminants,
we give summation formulas for Schur functions involving hyperdeterminants, and 
evaluate Toeplitz-type hyperdeterminants by using Jack function theory.

\medskip

\noindent
{\bf MSC-class}: primary 05E05; secondary 15A15 \\
{\bf Key words}:  Jacobi-Trudi formula, hyperdeterminant, Jack function,
 Schur function, Toeplitz determinant, pfaffian.
\end{abstract}


%
\section{Introduction}
%

The Schur function associated with a partition $\lambda$ of length $\le n$
has the determinantal expression
$$
s_{\lambda} = \det (h_{\lambda_i-i+j})_{1 \le i,j \le n},
$$
where $h_k$ is the complete symmetric function, i.e., the one-row Schur function $s_{(k)}$.
This formula is called the Jacobi-Trudi formula.

Jack functions are symmetric functions indexed by partitions,
i.e., Young diagrams.
They have one parameter $\alpha>0$ 
and include Schur functions as the special case $\alpha=1$.
In this paper, we consider the problem
of obtaining a Jacobi-Trudi type formula for Jack functions.
Lassalle and Schlosser \cite{Lassalle,LS} have recently 
derived such a formula for Macdonald functions.  
Their formula gives an expansion of any Macdonald function
with respect to one-row Macdonald functions or elementary symmetric functions.
Macdonald functions \cite{Mac} are generalizations of Jack functions, and so
the formula of Lassalle and Schlosser reduces to one for Jack functions as a limit case.
However, their formula is very complicated and is not expressed in a determinant-like form.
Kerov obtains a determinant expression for
Macdonald functions (and so Jack functions) of hook shapes \cite{Kerov}.
Also, the other determinantal expression is seen in \cite{Lapointe}.

We would like to obtain a Jacobi-Trudi formula for Jack functions
expressed in a determinant-like form.
In particular, the goal of this paper is to obtain the Jacobi-Trudi type formula
for Jack functions of {\it rectangular shapes} 
when the parameter $\alpha$ is either an integer $m$
or its inverse $1/m$.
By Jack functions of a rectangular shape 
we mean the Jack functions associated with a rectangular Young diagram.

In our formula, we employ a {\it hyperdeterminant}.
Cayley \cite{Cayley} defined some generalizations 
of the determinant to higher dimensional arrays.
Among them, we here deal with the following simple alternating sum 
$$
\hdet{2m}(A) := \frac{1}{n!} \sum_{\sigma_1, \dots, \sigma_{2m} \in \mf{S}_n}
\sgn(\sigma_1) \cdots \sgn(\sigma_{2m}) \prod_{i=1}^n A(\sigma_1(i), \dots, \sigma_{2m} (i))
$$
for an array $A=(A(i_1,\dots,i_{2m}))_{1 \le i_1, \dots, i_{2m} \le n}$.
We call it the hyperdeterminant of $A$.
This polynomial has been studied in \cite{Haukkanen, Barvinok,
LuqueThibon2003, LuqueThibon2004, Luque}, see also \cite{Sokolov}. 
Specifically,
a Cauchy-Binet type summation formula 
for hyperdeterminants is obtained (see e.g. \cite{Barvinok}).
Luque and Thibon \cite{LuqueThibon2003, LuqueThibon2004} 
studied the hyperdeterminant analogue of Hankel determinants,
which is closely related to Selberg's integral evaluation.
Furthermore, Haukkanen \cite{Haukkanen} and Luque \cite{Luque} studied
a hyperdeterminant analogue of GCD matrices and its extension to semilattices.

Our main result is stated as follows:
for the Jack $Q$-function $Q^{(1/m)}_{(L^n)}$ associated with the partition 
$(L,L,\dots,L)$ with $n$ components
and parameter $\alpha=1/m$, we obtain the following expression:
\begin{equation} \label{EqJackJT}
Q_{(L^n)}^{(1/m)} = \frac{n! \, (m!)^n}{(mn)!} \hdet{2m}
(g^{(1/m)}_{L+i_1+ \cdots +i_m-i_{m+1}-\cdots-i_{2m}})_{1 \le i_1,\dots,i_{2m} \le n},
\end{equation}
where $g^{(1/m)}_k$ are one-row Jack $Q$-functions.
If we substitute $m=1$ in expression \eqref{EqJackJT},
we obtain the Jacobi-Trudi formula for Schur functions of rectangular shapes
$s_{(L^n)}= \det(h_{L+i-j})_{1 \le i,j \le n}$.
Thus, our formula \eqref{EqJackJT} is 
the Jacobi-Trudi formula for Jack functions of rectangular shapes.

We give some applications of hyperdeterminants here.
In detail, this paper is organized as follows:

In \S 2, we develop the theory of hyperdeterminants and hyperpfaffians, 
where the hyperpfaffian is a pfaffian analogue of hyperdeterminants.
We give a Cauchy-Binet type integral formula involving hyperdeterminants,
which was previously essentially obtained in summation form.
However, the present integral form yields many applications.
We also define a hyperpfaffian
and give related integral formulas.
These integral formulas will be applied in later sections.
Note that our hyperpfaffian differs from that of Barvinok;
the explicit relation between these two hyperpfaffians 
 is given in Appendix \S 6.

In \S 3, we obtain
summation formulas for Schur functions by employing the formulas obtained in \S 2.
We express the summation
$$
\sum_{\lambda:\ell(\lambda) \le n} 
s_{\lambda}(\bx^{(1)}) s_{\lambda}(\bx^{(2)}) \cdots s_{\lambda}(\bx^{(k)}),
\qquad \text{where each $\bx^{(i)}$ is the set of $n$ variables,}
$$
as a hyperdeterminant if $k$ is even, or as a hyperpfaffian if $k$ is odd.
Our formulas in \S 3 are generalizations of Cauchy's determinant formula, 
Gessel's formula \cite{Gessel}, and Stembridge's formula \cite{Stembridge}.

In \S 4,
we compute Toeplitz hyperdeterminants.
A Toeplitz hyperdeterminant is a hyperdeterminant whose entries are given by
$d(i_1+i_2+ \cdots +i_m -i_{m+1} - \cdots -i_{2m})$,
where $d(k)$ are Fourier coefficients for a function defined on the unit circle.
We compute values of some Toeplitz hyperdeterminants explicitly
by employing the theory of Jack functions and the integral formulas obtained in \S 2.
 The problem of calculating a Toeplitz hyperdeterminant is then reduced to
the problem of evaluating a coefficient in the expansion of a certain symmetric  function
with respect to Jack functions.
Furthermore, we obtain a Toeplitz hyperdeterminant version of
the strong Szeg\"{o} limit theorem.
We derive the
 asymptotic behavior of a Toeplitz hyperdeterminant 
when the dimension $n$ goes to the infinity.

Finally, we prove our main result \eqref{EqJackJT} in \S 5,
employing the technique for Toeplitz hyperdeterminants developed in \S 4.

\section{Hyperdeterminants and hyperpfaffians}

In this section, we first give a Cauchy-Binet type integral formula
for hyperdeterminants.
 We define a new type of hyperpfaffian and obtain some integral formulas
making use of this concept.
These formulas will be applied in later sections.

We first state the definition of the hyperdeterminant again,
see also \cite{Cayley}.
Let $m$ and $n$ be positive integers.
For an array $A=(A(i_1,\dots,i_{2m}))_{1 \le i_1, \dots, i_{2m} \le n}$,
the hyperdeterminant of $A$ is defined by the alternating sum
$$
\hdet{2m}(A(i_1,\dots,i_{2m}))_{1 \le i_1, \dots, i_{2m} \le n}  \\
:=  \frac{1}{n!} \sum_{\sigma \in \mf{S}_n^{2m}}
\sgn(\sigma) 
\prod_{i=1}^n A(\sigma_1(i), \sigma_2(i), \dots, \sigma_{2m} (i)),
$$
where $\sgn(\sigma) =\sgn(\sigma_1) \cdots \sgn(\sigma_{2m})$ if 
$\sigma=(\sigma_1,\dots, \sigma_{2m}) \in \mf{S}_n^{2m}$.
We will sometimes write $(A(i_1,\dots,i_{2m}))_{1 \le i_1, \dots, i_{2m} \le n}$
as $(A(i_1,\dots,i_{2m}))_{[n]}$
for short.

The following proposition is a generalization of the Cauchy-Binet formula,
and was previously essentially obtained in  summation form,
see e.g. \cite{Barvinok}.
The case $m=1$ is the well-known Cauchy-Binet formula for determinants.

\begin{prop} \label{Prop:CauchyBinet}
Let $(X,\mu(\dd x))$ be a measure space 
and $\{\phi_{i,j}\}_{1 \le i \le 2m, 1 \le j \le n}$ 
functions on $X$.
Assume
$$
M(i_1,\dots, i_{2m}) := \int_X \phi_{1,i_1}(x) \phi_{2,i_2}(x) \cdots 
\phi_{2m,i_{2m}}(x) \, \mu(\dd x)
$$
is well-defined.
Then we have
$$
\frac{1}{n!} \int_{X^n} \prod_{i=1}^{2m} \det(\phi_{i,j}(x_k))_{1 \le j,k \le n} 
\cdot \prod_{j=1}^{n} \mu(\dd x_j)
=\hdet{2m}(M(i_1,\dots,i_{2m}))_{[n]}.
$$
\end{prop}

\begin{proof}
From definition of the hyperdeterminant, we see that
\begin{align*}
 \hdet{2m}(M(i_1,\dots,i_{2m}))_{[n]}  
=& \frac{1}{n!} \sum_{\sigma_1,\dots,\sigma_{2m} \in \mf{S}_n} 
\sgn(\sigma_1 \cdots \sigma_{2m}) \prod_{j=1}^n \( 
\int_X \prod_{i=1}^{2m}\phi_{i,\sigma_{i}(j)}(x) \, \mu(\dd x) \) \\
=& \frac{1}{n!} \int_{X^n}
\sum_{\sigma_1,\dots,\sigma_{2m} \in \mf{S}_n} 
\sgn(\sigma_1 \cdots \sigma_{2m}) \prod_{j=1}^n
\prod_{i=1}^{2m} \phi_{i,\sigma_{i}(j)}(x_j) \cdot \prod_{j=1}^{n} \mu(\dd x_j).
\end{align*}
Here, in the second equality, we switch the integral and the product.
The integrand on the right-hand side above is equal to
$$
\prod_{i=1}^{2m} \( \sum_{\sigma \in \mf{S}_n} \sgn(\sigma) \prod_{j=1}^{n}
\phi_{i,\sigma(j)}(x_j) \)
= \prod_{i=1}^{2m} \det (\phi_{i,k}(x_j))_{1 \le j,k \le n}, 
$$
and so we obtain the claim.
\end{proof}

\begin{remark}
We may obtain a permanent analogue of Proposition \ref{Prop:CauchyBinet} 
in a similar way.
For functions $\{ \phi_{i,j}\}_{1 \le i \le m, 1 \le j \le n}$,
it holds that
$$
\frac{1}{n!} \int_{X^n} \prod_{i=1}^{m} \hper{2}(\phi_{i,j}(x_k))_{1 \le j,k \le n} 
\cdot \prod_{j=1}^{n} \mu(\dd x_j)
=\hper{m}\( \int_X \prod_{k=1}^m 
\phi_{k,i_{k}}(x) \, \mu(\dd x) \)_{1 \le i_1,\dots,i_m  \le n}.
$$
Here $\hper{m}$ is a {\it hyper-permanent}
$$
\hper{m}(A(i_1,\dots,i_m))_{1 \le i_1,\dots,i_m \le n}
:= \frac{1}{n!} \sum_{\sigma_1,\dots,\sigma_m \in \mf{S}_n} \prod_{i=1}^n 
A(\sigma_1(i),\dots, \sigma_m(i)).
$$
\qed
\end{remark}

We next define a hyperpfaffian.
Let $B=(B(i_1,\dots,i_{2m}))_{[2n]}$ be an array
satisfying
\begin{equation} \label{EqAlter}
B(i_{\tau_1(1)}, i_{\tau_1(2)}, \dots, i_{\tau_m(2m-1)},i_{\tau_m(2m)})
= \sgn(\tau_1) \cdots \sgn(\tau_m) B(i_1,\dots,i_{2m})
\end{equation} 
for any $(\tau_1,\dots,\tau_m) \in (\mf{S}_2)^m$.
Here each $\tau_s \in \mf{S}_2$ permutes $2s-1$ with $2s$.
 We then define the  {\it hyperpfaffian} of $B$ by
\begin{equation} \label{Eq:DefHyperPfaffian}
\hpf{2m}(B) := \frac{1}{n!}\sum_{\sigma_1,\dots,\sigma_m \in \mf{E}_{2n}}
\sgn(\sigma_1 \cdots \sigma_{m}) \prod_{i=1}^n 
B(\sigma_1(2i-1),\sigma_1(2i), \dots, \sigma_m(2i-1),\sigma_m(2i)),
\end{equation}
where $\mf{E}_{2n}  :=
\{ \sigma \in \mf{S}_{2n} \ | \ \sigma(2i-1) < \sigma(2i) \ 
(1 \le i \le n) \}$.
When $m=1$, this expression is that of the ordinary pfaffian $\pf(B)$ of 
an alternating matrix $B=(B(i,j))_{1 \le i,j \le 2n}$.

\begin{remark}
Our hyperpfaffian is expressed as Barvinok's hyperpfaffian, see \S 6. \qed
\end{remark}

\begin{remark}
Let $\{\xi_i\}_{i \ge 1}$ be anti-commutative symbols, i.e., $\xi_i \xi_j= -\xi_j \xi_i$,
and let $\Lambda$ be the $\bC$-algebra generated by these $\xi_i$'s.
For a given array $A=(A(i_1,\dots,i_{2m}))_{[n]}$,
if we put 
$$
\eta=\sum_{1 \le k_1,\dots,k_{2m} \le n} A(k_1,\dots,k_{2m}) \, \xi_{k_1} \otimes \cdots \otimes
\xi_{k_{2m}} \in \Lambda^{\otimes 2m},
$$
then we have $\eta^n = n! \, \hdet{2m}(A) (\xi_1 \cdots \xi_{n})^{\otimes 2m}$.
Similarly, if we put
$$
\zeta= \sum_{1 \le k_1 < k_{2} \le 2n} \cdots \sum_{1 \le k_{2m-1} < k_{2m} \le 2n}
B(k_1,\dots,k_{2m}) (\xi_{k_1} \xi_{k_2}) \otimes \cdots \otimes (\xi_{k_{2m-1}} \xi_{k_{2m}}) 
\in \Lambda^{\otimes m}
$$
for a given tensor $B=(B(i_1,\dots,i_{2m}))_{[2n]}$
satisfying \eqref{EqAlter}, then $\zeta^n = n! \, \hpf{2m}(B) (\xi_1 \cdots \xi_{2n})^{\otimes m}$.
\qed 
\end{remark}

The following proposition describes the relationship between hyperpfaffians and hyperdeterminants. It will be shown that
any hyperdeterminant is a special case of a hyperpfaffian.

\begin{prop}
Let $A=(A(i_1,\dots, i_{2m}))_{[n]}$ be any array.
Define the array $B=(B(i_1,\dots,i_{2m}))_{[2n]}$ as follows.
If $i_{2s-1}$ is odd and $i_{2s}$ is even for all $1 \le s \le m$,
then 
$B(i_1,\dots,i_{2m})= A(p_1,q_1, \dots, p_m,q_m)$, 
where $i_{2s-1}=2p_s-1$ and $i_{2s}=2q_s$ for $1 \le s \le m$.
If $i_{2s-1}$ is even and $i_{2s}$ is odd for all $1 \le s \le m$,
then 
$B(i_1,\dots,i_{2m})$ is defined by \eqref{EqAlter}; otherwise
 $B(i_1,\dots,i_{2m})=0$. 
 We then have
$\hpf{2m}(B)= \hdet{2m}(A)$.
\end{prop}

\begin{proof}
From the alternating property \eqref{EqAlter} and
the definition of a hyperpfaffian, 
we can express a hyperpfaffian as alternating sums on $\mf{S}_n^m$;
$$
\hpf{2m}(B)= \frac{1}{2^{nm} \, n!}\sum_{\sigma_1,\dots,\sigma_m \in \mf{S}_{2n}}
\sgn(\sigma_1 \cdots \sigma_{m}) \prod_{i=1}^n 
B(\sigma_1(2i-1),\sigma_1(2i), \dots, \sigma_m(2i-1),\sigma_m(2i)).
$$
Because of the definition of $B(i_1,\dots,i_{2m})$,
each term vanishes
if there exist  $1 \le s \le m$ and $1 \le i \le n$ such that
${\displaystyle \sigma_s(2i-1) \equiv \sigma_s(2i) \pmod{2}}$.
We may therefore assume that on the sum of the above equality, 
only one of $\sigma_s(2i-1)$ and $\sigma_s(2i)$ is odd, say $2r_{2s-1}-1$, and 
another is even, say $2r_{2s}$ for any $1 \le s \le m$.
If $\sigma_{s}(2i-1)=2r_{2s}$ (and therefore 
${\displaystyle \sigma_{s}(2i)=2r_{2s-1}-1}$),
we replace $B(\dots, 2r_{2s}, 2r_{2s-1}-1, \dots)$ with
$-B(\dots, 2r_{2s-1}-1, 2r_{2s}, \dots)$. 
Then, for each ${\displaystyle 1 \le s \le m}$,
the sequences $(\sigma_s(1), \sigma_s(3), \dots, \sigma_s(2n-1))$ 
and $(\sigma_s(2), \sigma_s(4), \dots, \sigma_s(2n))$
are permutations of ${\displaystyle 1, 3, \dots, 2n-1}$ 
and $2, 4, \dots, 2n$, respectively. 
Hence we have the expression
$$
\hpf{2m}(B)= \frac{1}{2^{mn} n!}\sum_{\tau_1,\dots,\tau_{2m} \in \mf{S}_{n}}
\sgn(\sigma_1 \cdots \sigma_{m}) \prod_{i=1}^n 2^m
B(2 \tau_1(i)-1, 2 \tau_2(i), \dots, 2 \tau_{2m-1}(i)-1, 2 \tau_{2m}(i)),
$$
where $\sigma_s$ denotes the permutation
$$
\begin{pmatrix} 1 & 2& \cdots & 2n-1 &2n \\ 
2\tau_{2s-1}(1)-1 & 2\tau_{2s}(1) & \cdots & 2\tau_{2s-1}(n)-1 & 2\tau_{2s}(n)
\end{pmatrix} \in \mf{S}_{2n}.
$$
Since $\sgn(\sigma_s)= \sgn(\tau_{2s-1}) \sgn(\tau_{2s})$, we see that
$$
\hpf{2m}(B)=
\frac{1}{n!}\sum_{\tau_1,\dots,\tau_{2m} \in \mf{S}_{n}}
\sgn(\tau_1 \cdots \tau_{2m}) \prod_{i=1}^n 
A(\tau_1(i), \dots, \tau_{2m}(i))= \hdet{2m}(A).
$$
\end{proof}

We now give a hyperpfaffian analogue of Proposition \ref{Prop:CauchyBinet}.
The case $m=1$ is well known as the de Bruijn formula \cite{Bruijn}.
We will make use of this proposition in the proofs of Theorem \ref{Thm:OddSumSchur} and Theorem \ref{Thm:SumSchurGessel} below.

\begin{prop} \label{Prop:IntegralFormulaPf}
Let $m$ be an odd positive integer.
Let $\{\psi_{i,j}\}_{1 \le i \le m, 1 \le j \le 2n}$ be functions 
on a measure space $(X, \mu (\dd x))$ and $\epsilon$ be a function on $X \times X$
such that $\epsilon(y,x) = -\epsilon(x,y)$.
Suppose that
$$
Q(i_1,\dots,i_{2m})= \frac{1}{2} \int_{X^2} \epsilon(x,y)
\prod_{s=1}^m \det \begin{pmatrix} \psi_{s,i_{2s-1}}(x) & \psi_{s,i_{2s-1}}(y) \\
\psi_{s,i_{2s}}(x) & \psi_{s,i_{2s}}(y) \end{pmatrix}
\, \mu(\dd x) \mu(\dd y)
$$ 
is well-defined.
Then we have
$$
\frac{1}{(2n)!} \int_{X^{2n}} \pf (\epsilon(x_i,x_j))_{1 \le i,j \le 2n}
\prod_{s=1}^m \det(\psi_{s,j}(x_k))_{1 \le j,k \le 2n} \,
\prod_{j=1}^{2n} \mu(\dd x_j) 
= \hpf{2m}(Q(i_1,\dots,i_{2m}))_{[2n]}.
$$
\end{prop}

\begin{proof}
A straightforward calculation gives 
\begin{align*}
& \hpf{2m}(Q(i_1,\dots,i_{2m}))_{[2n]} \\ 
=&
\frac{1}{n!} \sum_{\sigma_1, \dots, \sigma_{m} \in \mf{E}_{2n}}
\sgn(\sigma_1 \cdots \sigma_m)
\prod_{j=1}^n \bigg\{
\frac{1}{2} \int_{X^2} \epsilon(x,y) \\
& \qquad \times 
\prod_{s=1}^m \det \begin{pmatrix} \psi_{s,\sigma_s(2j-1)}(x) & \psi_{s,\sigma_s(2j-1)}(y) \\
\psi_{s,\sigma_s(2j)}(x) & \psi_{s,\sigma_s(2j)}(y) \end{pmatrix}
\, \mu(\dd x) \mu(\dd y) \bigg\} \\
=& \frac{1}{n! \, 2^n} \int_{X^{2n}} \prod_{j=1}^n \epsilon(x_{2j-1},x_{2j}) \\
& \quad \times \prod_{s=1}^{m}
\left\{ \sum_{\sigma \in \mf{E}_{2n}} \sgn(\sigma) \prod_{j=1}^n
\det \begin{pmatrix} \psi_{s,\sigma(2j-1)}(x_{2j-1}) & \psi_{s,\sigma(2j-1)}(x_{2j}) \\
\psi_{s,\sigma(2j)}(x_{2j-1}) & \psi_{s,\sigma(2j)}(x_{2j}) \end{pmatrix} \right\}
\, \prod_{j=1}^{2n} \mu(\dd x_j).
\end{align*}
Since the well-known expansion (see e.g. Lemma 4 in \cite{IW})
$$
\det (a_{i,j})_{1 \le i,j \le {2n}} = 
\sum_{\sigma \in \mf{E}_{2n}} \sgn(\sigma) \prod_{j=1}^n
\det \begin{pmatrix} a_{\sigma(2j-1), 2j-1} & a_{\sigma(2j-1), 2j} \\
a_{\sigma(2j), 2j-1} & a_{\sigma(2j), 2j} \end{pmatrix},
$$
we have
\begin{equation} \label{EqPfQ1}
\hpf{2m}(Q(i_1,\dots,i_{2m}))_{[2n]} =
\frac{1}{n! \, 2^n} \int_{X^{2n}} \prod_{j=1}^n \epsilon(x_{2j-1},x_{2j}) \cdot \prod_{s=1}^{m}
\det (\psi_{s,j}(x_k))_{1 \le j,k \le 2n} \, \prod_{j=1}^{2n} \mu(\dd x_j).
\end{equation}

On the other hand, by expanding the pfaffian, we see that
\begin{align*}
& \frac{1}{(2n)!} \int_{X^{2n}} \pf (\epsilon(x_i,x_j))_{1 \le i,j \le 2n}
\prod_{s=1}^m \det(\psi_{s,j}(x_k))_{1 \le j,k \le 2n} \,
\prod_{j=1}^{2n} \mu(\dd x_j) \\
=& \frac{1}{(2n)! \, n!} \sum_{\sigma \in \mf{E}_{2n}} \int_{X^{2n}}
\sgn(\sigma) \prod_{j=1}^n \epsilon(x_{\sigma(2j-1)}, x_{\sigma(2j)}) \cdot 
\prod_{s=1}^m \det(\psi_{s,j}(x_k))_{1 \le j,k \le 2n} \,
\prod_{j=1}^{2n} \mu(\dd x_j).
\end{align*}
Since $m$ is odd, by permuting columns of 
each $\det(\psi_{s,j}(x_{k}))_{1 \le j,k \le 2n}$,
we have
\begin{align*}
=& \frac{1}{(2n)! \, n!} \sum_{\sigma \in \mf{E}_{2n}} \int_{X^{2n}}
\prod_{j=1}^n \epsilon(x_{\sigma(2j-1)}, x_{\sigma(2j)}) \cdot 
\prod_{s=1}^m \det(\psi_{s,j}(x_{\sigma(k)}))_{1 \le j,k \le 2n} \,
\prod_{j=1}^{2n} \mu(\dd x_j) \\
=& \frac{1}{n! \, 2^n} \int_{X^{2n}}
\prod_{j=1}^n \epsilon(x_{2j-1}, x_{2j}) \cdot 
\prod_{s=1}^m \det(\psi_{s,j}(x_{k}))_{1 \le j,k \le 2n} \,
\prod_{j=1}^{2n} \mu(\dd x_j).
\end{align*}
Here we use the fact the cardinality of the set $\mf{E}_{2n}$ is $(2n)!/2^n$.
Combining the above expression with equation \eqref{EqPfQ1}, we obtain the claim.
\end{proof}

We also obtain another integral formula in a similar way to the previous proposition.
The following proposition will be used in the proof of Theorem \ref{Thm:ToeplitzPfaffian}.

\begin{prop} \label{Prop:IntegralFormulaPf2}
Let $\{\psi_{i,j}\}_{1 \le i \le 2m, 1 \le j \le 2n}$ be functions on a measure space 
$(X,\mu(\dd x))$.
Suppose that
$$
R(i_1,\dots,i_{2m})= \int_X \prod_{s=1}^m \det
\begin{pmatrix} \psi_{2s-1,i_{2s-1}}(x) & \psi_{2s-1,i_{2s}}(x) \\
\psi_{2s,i_{2s-1}}(x) & \psi_{2s,i_{2s}}(x) \end{pmatrix}
\mu(\dd x)
$$
is well-defined.
Then we have
$$
\frac{1}{n!} \int_{X^n} \prod_{s=1}^m 
\det ( \psi_{2s-1,j}(x_k) \ | \ \psi_{2s,j}(x_k))_{1 \le j \le 2n, 1\le  k \le n} \
\prod_{j=1}^n \mu(\dd x_j) 
=\hpf{2m}(R(i_1,\dots,i_{2m}))_{[2n]}.
$$
Here $\det(a_{j,k} \ | \ b_{j,k})_{1\le j \le 2n, 1\le  k \le n}$ denotes
the determinant of the $2n$ by $2n$ matrix whose $j$-th row is given by 
$(a_{j,1} \ b_{j,1} \ a_{j,2} \ b_{j,2} \ \dots \ a_{j,n} \ b_{j,n})$.   
\qed
\end{prop}

%
\section{Summation formulas for Schur functions}
%

From the propositions obtained in the previous section,
we obtain several summation formulas for Schur functions.
We express the summation of the product of Schur functions
as a hyperdeterminant or hyperpfaffian.

Let $\bx=(\bx_1,\dots, \bx_n)$ be a sequence of $n$ variables. 
For a sequence $\alpha=(\alpha_1,\dots,\alpha_n)$ of non-negative integers,
we introduce $a_{\alpha}(\bx) = \det(\bx_i^{\alpha_j})_{1 \le i,j \le n}$ as in \cite{Mac}.
In particular, we have 
${\displaystyle a_{\delta}(\bx)= \det(\bx_i^{n-j})_{1 \le i,j \le n}=V(\bx)}$, 
where $V(\bx)=\prod_{1 \le i < j \le n} (\bx_i-\bx_j)$ stands for 
the Vandermonde product 
and  $\delta=\delta_n=(n-1,n-2,\dots, 1,0)$.
For a partition $\lambda=(\lambda_1 \ge \dots \ge \lambda_n)$ of length $\le n$,
the Schur function $s_{\lambda}$ corresponding to $\lambda$
is defined by $s_{\lambda}(\bx) = 
a_{\lambda+\delta}(\bx) / V(\bx)$, where 
$\lambda+\delta=(\lambda_1+n-1, \lambda_2+n-2, \dots, \lambda_n)$.

The Schur functions have the following well-known summation formulas
\cite[\S I-4, Ex. 6]{Mac}, \cite{IOW}:
\begin{align}
\sum_{\lambda:\ell(\lambda) \le n} s_{\lambda}(\bx) s_{\lambda}(\by) 
=& \frac{1}{V(\bx) V(\by)} \cdot 
\det \( \frac{1}{1-\bx_i \by_j} \)_{1 \le i,j \le n}, 
\qquad \by=(\by_1, \dots, \by_n), \label{EqCauchySumDet} \\
\sum_{\lambda:\ell(\lambda) \le n} s_{\lambda}(\bx) 
=& \frac{1}{V(\bx)}  \cdot
\pf \( \frac{\bx_i- \bx_j}{(1-\bx_i) (1-\bx_j)(1-\bx_i \bx_j) } \)_{1 \le i,j \le n}. 
\label{EqLittlewoodSumPf}
\end{align}
Here  we assume $n$ is even in expression \eqref{EqLittlewoodSumPf}.
Note that the determinant on the right-hand side of equation \eqref{EqCauchySumDet}
is called Cauchy's determinant.
We extend expressions \eqref{EqCauchySumDet} and 
\eqref{EqLittlewoodSumPf}
to higher degrees.

\begin{thm} \label{Thm:CauchyDetFinite}
Let $\bx^{(i)}=(\bx^{(i)}_1, \bx^{(i)}_2,\dots, \bx^{(i)}_n)$ 
be a sequence of variables for each $1 \le i \le 2m$.
For any positive integer $N$,
we have
$$
\sum
s_{\lambda}(\bx^{(1)}) s_{\lambda}(\bx^{(2)}) \cdots s_{\lambda}(\bx^{(2m)})
=\prod_{i=1}^{2m} \frac{1}{V(\bx^{(i)})} \cdot
\hdet{2m}\( 
\frac{1-(\bx_{i_1}^{(1)} \cdots \bx_{i_{2m}}^{(2m)})^{n+N}}%
{1-\bx_{i_1}^{(1)} \cdots \bx_{i_{2m}}^{(2m)}} \)_{[n]},
$$ 
where the sum is over all partitions $\lambda=(\lambda_1, \dots, \lambda_n)$
of length $\ell(\lambda) \le n$ and
of largest part $\lambda_1 \le N$.
\end{thm}

\begin{proof}
Let $X=\{0,1,\dots, n+N-1\}$ and $\phi_{i,j}(k)= (\bx_{j}^{(i)})^k$,
and apply Proposition \ref{Prop:CauchyBinet}.  
Here the integrals in Proposition \ref{Prop:CauchyBinet}
are regarded as summations over $X$.  
Then, since
$$
\sum_{k=0}^{n+N-1} \phi_{1,i_1}(k) \phi_{2,i_2}(k) \cdots \phi_{2m,i_{2m}}(k)
= \sum_{k=0}^{n+N-1} (\bx_{i_1}^{(1)} \cdots \bx_{i_{2m}}^{(2m)})^k
= \frac{1-(\bx_{i_1}^{(1)} \cdots \bx_{i_{2m}}^{(2m)})^{n+N}}%
{1-\bx_{i_1}^{(1)} \cdots \bx_{i_{2m}}^{(2m)}},
$$
we have
$$
\hdet{2m}\( 
\frac{1-(\bx_{i_1}^{(1)} \cdots \bx_{i_{2m}}^{(2m)})^{n+N}}%
{1-\bx_{i_1}^{(1)} \cdots \bx_{i_{2m}}^{(2m)}} \)_{[n]}
= \sum_{n+N-1 \ge j_1 > \dots > j_n \ge 0} \prod_{i=1}^{2m} 
\det( (\bx_{p}^{(i)})^{j_q})_{1 \le p,q \le n}.
$$
Replacing each sequence $(j_1, \dots, j_n)$ with a partition $\lambda$ by 
$\lambda_q +n-q = j_q$, 
we see that 
$$
\hdet{2m}\( 
\frac{1-(\bx_{i_1}^{(1)} \cdots \bx_{i_{2m}}^{(2m)})^{n+N}}%
{1-\bx_{i_1}^{(1)} \cdots \bx_{i_{2m}}^{(2m)}} \)_{[n]}
=\sum_{\begin{subarray}{c} \lambda: \ell(\lambda) \le n, \\ \lambda_1 \le N \end{subarray}}
a_{\lambda+\delta}(\bx^{(1)}) a_{\lambda+\delta}(\bx^{(2)})
\cdots a_{\lambda+\delta}(\bx^{(2m)}) 
$$
and the claim follows.
\end{proof}

\begin{cor} \label{Cor:CauchyIdentityHyper}
Let $\bx^{(i)}$ be as in Theorem \ref{Thm:CauchyDetFinite}. Then
\begin{equation} \label{Eq:CauchyIdentityHyper}
\sum_{\lambda: \ell(\lambda) \le n}
s_{\lambda}(\bx^{(1)}) s_{\lambda}(\bx^{(2)}) \cdots s_{\lambda}(\bx^{(2m)})
= \prod_{i=1}^{2m} \frac{1}{V(\bx^{(i)})} \cdot 
\hdet{2m}\( 
\frac{1}%
{1-\bx_{i_1}^{(1)} \cdots \bx_{i_{2m}}^{(2m)}} \)_{[n]}.
\end{equation}
\end{cor}

\begin{proof}
We may assume that each variable $\bx_j^{(i)}$ belongs to 
the open unit disc $\{z \in \bC \ | \ |z|<1\}$.
We obtain the claim by taking the limit 
as $N \to \infty$
in Theorem \ref{Thm:CauchyDetFinite}.
\end{proof}

One may see expression \eqref{Eq:CauchyIdentityHyper} as 
a simple multi-version of expression \eqref{EqCauchySumDet}.
An odd-product analogue of expression \eqref{Eq:CauchyIdentityHyper} is given as follows:

\begin{thm} \label{Thm:OddSumSchur}
Let $m$ be an odd positive number and 
$\bx^{(i)}=(\bx^{(i)}_1, \bx^{(i)}_2,\dots, \bx^{(i)}_{2n})$ 
 the sequence of variables for each $1 \le i \le m$.
Then we have
$$
 \sum_{\lambda: \ell(\lambda) \le 2n}
s_{\lambda}(\bx^{(1)}) s_{\lambda}(\bx^{(2)}) \cdots s_{\lambda}(\bx^{(m)}) 
= \prod_{i=1}^{m} \frac{1}{V(\bx^{(i)})} \cdot
\hpf{2m} \( \frac{\sum_{p=0}^\infty \prod_{s=1}^m \{
(\bx_{i_{2s-1}}^{(s)})^{p+1} - (\bx_{i_{2s}}^{(s)})^{p+1}\} 
 }{
1-\bx_{i_1}^{(1)} \bx_{i_2}^{(1)} \cdots \bx_{i_{2m-1}}^{(m)} \bx_{i_{2m}}^{(m)}} \)_{[2n]}.
$$
\end{thm}

\begin{proof}
Let $X=\{0,1,2,\dots\}$ and $\psi_{i,j}(k)= (\bx_j^{(i)})^k$.
Let $\epsilon$ be the alternating function defined by $\epsilon(k,l)=1$ for $k>l$ and 
apply Proposition \ref{Prop:IntegralFormulaPf}. 
Then we have
\begin{align*}
& Q(i_1,\dots,i_{2m}) = \sum_{k>l \ge 0} \prod_{s=1}^m 
\{ (\bx_{i_{2s-1}}^{(s)})^k (\bx_{i_{2s}}^{(s)})^l 
-(\bx_{i_{2s-1}}^{(s)})^l (\bx_{i_{2s}}^{(s)})^k \} \\
=& \sum_{p,l \ge 0} 
\prod_{s=1}^m ( \bx_{i_{2s-1}}^{(s)} \bx_{i_{2s}}^{(s)})^l 
\{ (\bx_{i_{2s-1}}^{(s)})^{p+1} -(\bx_{i_{2s}}^{(s)})^{p+1} \} 
= \frac{\sum_{p=0}^\infty \prod_{s=1}^m \{
(\bx_{i_{2s-1}}^{(s)})^{p+1} - (\bx_{i_{2s}}^{(s)})^{p+1}\} }{
1-\bx_{i_1}^{(1)} \bx_{i_2}^{(1)} \cdots \bx_{i_{2m-1}}^{(m)} \bx_{i_{2m}}^{(m)}}.
\end{align*}
 Proposition \ref{Prop:IntegralFormulaPf} therefore implies that
\begin{align*}
\hpf{2m}(Q(i_1,\dots,i_{2m}))_{[2n]}=& \sum_{k_1 > k_2 > \cdots >k_{2n} \ge 0}
\pf ( \epsilon(k_i,k_j))_{1 \le i,j \le 2n} 
\prod_{s=1}^m \det ( (\bx_{p}^{(s)})^{k_q})_{1 \le p,q \le 2n} \\
=& \sum_{\lambda_1 \ge \cdots \ge \lambda_{2n} \ge 0} 
\prod_{s=1}^m \det ( (\bx_{p}^{(s)})^{\lambda_q+2n-q})_{1 \le p,q \le 2n}.
\end{align*}
Here we have replaced each $(k_1,\dots,k_{2N})$ 
with a partition $(\lambda_1,\dots,\lambda_{2N})$
by $k_{q}=\lambda_q+2N-q$. 
Thus, the theorem follows.
\end{proof}

Recall the Jacobi-Trudi formula for Schur functions
\begin{equation} \label{Eq:JacobiTrudi}
s_{\lambda}= \det (h_{\lambda_i-i+j})_{1 \le i,j \le n}
\end{equation}
for any partition $\lambda$ of length $\le n$.
Here $h_k$ is the complete symmetric function 
$$
h_k(\bx) = \sum_{\begin{subarray}{c} k_1,k_2, \dots \ge 0, \\
k_1+k_2 +\cdots=k \end{subarray}} \bx_1^{k_1} \bx_2^{k_2} \cdots
$$
in countably many variables $\bx=(\bx_1,\bx_2,\dots)$ if $k \ge 0$,
or $h_k=0$ otherwise.
From equation \eqref{Eq:JacobiTrudi} we obtain another summation formula 
for Schur functions.

\begin{thm} \label{Thm:SumSchurGessel}
Let  $\bx^{(i)}=(\bx^{(i)}_1, \bx^{(i)}_2,\dots)$ for each $i \ge 1$.
Then, for each positive integer $m$, we have
\begin{equation} \label{Eq:GesselHyper}
\sum_{\lambda: \ell(\lambda) \le n} 
s_{\lambda}(\bx^{(1)})  \cdots s_{\lambda}(\bx^{(2m)})
= \hdet{2m} \( \sum_{ k \ge 0} h_{k-i_1}(\bx^{(1)}) h_{k-i_2}(\bx^{(2)})
\cdots h_{k-i_{2m}}(\bx^{(2m)})\)_{[n]}.
\end{equation}
Furthermore, for an odd positive integer $m$, 
\begin{align} 
& \sum_{\lambda: \ell(\lambda) \le 2n} 
s_{\lambda}(\bx^{(1)})  \cdots s_{\lambda}(\bx^{(m)}) \label{Eq:StembridgeHyper}  \\
=& \hpf{2m} \( \sum_{k,l \ge 0} \prod_{s=1}^m
 \det \begin{pmatrix} 
h_{k+l+1-i_{2s-1}}(\bx^{(s)}) & h_{l-i_{2s-1}}(\bx^{(s)}) \\
h_{k+l+1-i_{2s}}(\bx^{(s)}) & h_{l-i_{2s}}(\bx^{(s)}) 
\end{pmatrix}
\)_{[2n]}. \notag
\end{align}
\end{thm}

\begin{proof}
Apply Proposition \ref{Prop:CauchyBinet} to $X=\{0,1,2,\dots\}$ and 
$\phi_{i,j}(k) = h_{k-j+1}(\bx^{(i)})$.
Then we have
\begin{align*}
\sum_{k_1 > \cdots >k_n \ge 0} \prod_{i=1}^{2m} 
\det (h_{k_q-p+1}(\bx^{(i)}))_{1 \le p,q \le n} =& 
\sum_{\lambda_1 \ge  \cdots \ge \lambda_n \ge 0} \prod_{i=1}^{2m} 
\det (h_{\lambda_q+n-q-(n-p+1)+1}(\bx^{(i)}))_{1 \le p,q \le n}  \\
=& \sum_{\lambda: \ell(\lambda) \le n} 
\prod_{i=1}^{2m} \det (h_{\lambda_q-q+p}(\bx^{(i)}))_{1 \le p,q \le n}.
\end{align*}
Here in the first step we have replaced each $k_q$ with $\lambda_q+n-q$ by a partition $\lambda$
and changed the order of rows of determinants $\det (h_{k_q-p+1}(\bx^{(i)}))_{1 \le p,q \le n}$.
Hence
the claim \eqref{Eq:GesselHyper} follows from expression \eqref{Eq:JacobiTrudi}.
The second formula \eqref{Eq:StembridgeHyper} 
follows by applying Proposition \ref{Prop:IntegralFormulaPf} 
to $X=\{0,1,\dots\}$, $\epsilon(k,l)=1$ for $k>l$, and $\psi_{i,j}(k)=h_{k-j+1}(\bx^{(i)})$
in a similar way.
\end{proof}

The cases $m=1$ in expressions \eqref{Eq:GesselHyper} and \eqref{Eq:StembridgeHyper} have previously been obtained in 
\cite{Gessel} and \cite{Stembridge} respectively.

%
\section{Toeplitz hyperdeterminants}
%

In this section,
we consider a class of hyperdeterminants called Toeplitz hyperdeterminants,
and we evaluate them
by employing the theory of Jack polynomials.
Furthermore, we obtain the strong Szeg\"{o} limit formula for Toeplitz hyperdeterminants.

\subsection{Heine-Szeg\"{o} formula for Toeplitz hyperdeterminants}

Let $f(z)$ be a complex-valued function on the unit circle $\bT=\{z \in \bC \ | \ |z|=1 \}$
whose Fourier expansion is given by $f(z) =\sum_{k \in \bZ} d(k) z^k$.
Then the hyperdeterminant
$$
D^{[2m]}_{n}(f) = \hdet{2m}(d(i_1+ \cdots +i_m - i_{m+1} - \cdots -i_{2m}))_{[n]}
$$
is called the {\it Toeplitz hyperdeterminant} of $f$,
see \cite{LuqueThibon2003}.

\begin{thm} \label{Thm:HSF}
For a function $f \in L^1 (\bT)$, we have
$$
D^{[2m]}_{n}(f) = \frac{1}{n!} \int_{\bT^n} f(z_1) f(z_2) \cdots f(z_n) 
|V(z_1,\dots,z_n)|^{2m} \, \dd z_1 \cdots \dd z_n,
$$
where $\dd z_j$ is the Haar measure on $\bT$ normalized by $\int_{\bT}\dd z_j=1$. 
\end{thm}

\begin{proof}
Apply Proposition \ref{Prop:CauchyBinet} to functions 
$\{\phi_{i,j}\}_{1 \le i \le 2m, 1 \le j \le n}$ on the measure space $(\bT, \dd z)$,
where
$$
\phi_{i,j}(z) = 
\begin{cases} 
f(z) z^{n-j} & \text{for $i=1$},\\
z^{n-j} & \text{for $2 \le i \le m$},\\
z^{j-n} & \text{for $m+1 \le i \le 2m$}.
\end{cases}
$$
Then, since
$\int_{\bT} \phi_{1,i_1} (z) \cdots \phi_{2m,i_{2m}} (z) \, \dd z
= d(i_1+ \cdots +i_{m} - i_{m+1} - \cdots - i_{2m})$,
we have
$$
D_{n}^{[2m]}(f) = \frac{1}{n!} \int_{\bT^n} 
f(z_1) \cdots f(z_n) \{ \det (z_k^{n-j})_{1 \le j,k \le n}
\det (z_k^{j-n})_{1 \le j,k \le n} \}^m \,
\dd z_1 \cdots \dd z_n,
$$
and upon using $\det (z_k^{n-j})_{1 \le j,k \le n}
\det (z_k^{j-n})_{1 \le j,k \le n} =|V(z_1,\dots, z_n)|^2$,
we have the claim.
\end{proof}

The case $m=1$ in
this theorem is simply 
the Heine-Szeg\"{o} formula for a Toeplitz determinant,
see e.g. \cite{BD}.

We can express $D^{[4m]}_{n}(f)$ by a hyperpfaffian.

\begin{thm} \label{Thm:ToeplitzPfaffian}
Toeplitz hyperdeterminants $D^{[4m]}_n(f)$ can be a hyperpfaffian $\hpf{2m}$.
$$
D^{[4m]}_{n}(f) = 
\hpf{2m} \( \prod_{s=1}^m (i_{2s}-i_{2s-1}) \cdot d((2n+1)m-\sum_{k=1}^{2m}i_k) \)_{[2n]}.
$$
\end{thm}

\begin{proof}
Apply Proposition \ref{Prop:IntegralFormulaPf2} to functions 
$\{\psi_{i,j}\}_{1 \le i \le 2m, 1 \le j \le n}$ on the measure space $(\bT, \dd z)$,
where
$$
\psi_{i,j}(z) = 
\begin{cases} 
f(z) z^{j-n-\frac{1}{2}} & \text{if $i=1$}, \\
z^{j-n-\frac{1}{2}} & \text{if $i$ is odd and $i>1$}, \\
(j-n-\frac{1}{2} )z^{j-n-\frac{1}{2}} & \text{if $i$ is even}. \end{cases}
$$
Then we have
$$
\hpf{2m}(R)= \frac{1}{n!} \int_{\bT^n} f(z_1) \cdots f(z_n) 
|V(z_1,\dots, z_n)|^{4m} \ \dd z_1 \cdots \dd z_n
=D^{[4m]}_{n}(f)
$$
upon using the formula (see e.g. \cite[Chapter 11]{Mehta})
$$
\det\( z^{j-n-\frac{1}{2}}_k \ \Bigm| \ (j-n-\frac{1}{2}) z^{j-n-\frac{1}{2}}_k
\)_{1 \le j \le 2n, 1 \le  k \le n}
= \prod_{1 \le j< k \le n} |z_j-z_k|^{4}.
$$
Here the entry of $R$ is calculated as
\begin{align*}
& R(i_1,\dots,i_{2m}) = \int_{\bT} f(z) \prod_{s=1}^m
\{ (i_{2s}-n-\frac{1}{2})- (i_{2s-1}-n-\frac{1}{2})\} z^{i_{2s-1}+i_{2s}-2n-1} \ \dd z \\
=& \prod_{s=1}^m (i_{2s}-i_{2s-1})  \int_{\bT} f(z) z^{i_1+\cdots+i_{2m}-m(2n+1)} \dd z
= \prod_{s=1}^m (i_{2s}-i_{2s-1}) \cdot d(m(2n+1)-i_1- \cdots- i_{2m})
\end{align*}
and so we obtain the claim.
\end{proof}

In particular, it holds that $D^{[4]}_{n}(f) = \pf \( (j-i) d(2n+1-i-j) \)_{1 \le i,j \le 2n}$.
For example,
$$
D^{[4]}_1(f)= d(0), \qquad D^{[4]}_2(f)= d(2)d(-2)-4d(1)d(-1)+3d(0)^2, \dots.
$$

\subsection{Basic properties of Jack functions}

To compute Toeplitz hyperdeterminants,
we employ Jack polynomials.
We recall the basic properties of Jack functions, see \cite[\S VI-10]{Mac} for details.
Let $\alpha>0$ and 
let $\Lambda(\alpha)$ be the $\bQ(\alpha)$-algebra of symmetric functions
with variables $\bx=(\bx_1,\bx_2,\dots,)$.
Let $p_k$ be the power-sum polynomial $p_k(\bx)= \bx_1^k+\bx_2^k+ \cdots $
and put $p_{\lambda}=p_{\lambda_1} p_{\lambda_2}\cdots $ for a partition 
$\lambda=(\lambda_1,\lambda_2,\dots)$.
Define the scalar product on $\Lambda(\alpha)$ by 
$$
\langle p_{\lambda}, p_{\mu} \rangle_{\alpha} = \delta_{\lambda,\mu} z_{\lambda}
\alpha^{\ell(\lambda)}
$$
for partitions $\lambda$ and  $\mu$.
Here $\delta_{\lambda,\mu}$ is Kronecker's delta and
$z_{\lambda}=\prod_{k \ge 1} k^{m_k} m_k!$ where 
${\displaystyle m_k=\#\{i \ge 1 \ | \ \lambda_i=k\}}$.
We denote a rectangular-shape partition $(k,k, \dots,k)$ with $n$ components by $(k^n)$.

Let $m_\lambda$ be the monomial symmetric function.
Then Jack $P$-functions $P_{\lambda}^{(\alpha)}$ are characterized as homogeneous 
symmetric functions such that
$$
P_{\lambda}^{(\alpha)} =m_{\lambda} + \sum_{\mu < \lambda} u_{\lambda\mu} m_{\mu}
\quad \text{with $u_{\lambda\mu} \in \bQ(\alpha)$},
\qquad \langle P^{(\alpha)}_{\lambda}, P^{(\alpha)}_{\mu} \rangle_{\alpha} =0
\quad \text{if $\lambda \not=\mu$},
$$
where ``$<$'' is the dominance ordering.
Put 
$$
c_{\lambda}(\alpha)= \prod_{(i,j) \in \lambda} (\alpha(\lambda_i-j) +\lambda_j'-i+1)
\qquad \text{and} \qquad
c'_{\lambda}(\alpha)= \prod_{(i,j) \in \lambda} (\alpha(\lambda_i-j+1) +\lambda_j'-i),
$$
where $(i,j)$ run over all squares in the Young diagram associated with $\lambda$
and $\lambda'=(\lambda_1', \lambda_2',\dots)$ is the conjugate partition of $\lambda$.
The set of functions $\{P_{\lambda}^{(\alpha)} \ | \ \text{$\lambda$ are partitions} \}$
is an orthogonal basis of $\Lambda(\alpha)$.
Defining Jack $Q$-functions by $Q_{\lambda}^{(\alpha)} = c_{\lambda}(\alpha)
c'_{\lambda}(\alpha)^{-1} P_{\lambda}^{(\alpha)}$,
the set  $\{Q_{\lambda}^{(\alpha)} \ | \ \text{$\lambda$ are partitions} \}$ is its dual basis,
i.e., $\langle P_{\lambda}^{(\alpha)}, Q_{\mu}^{(\alpha)} \rangle_{\alpha}=\delta_{\lambda,\mu}$.

We sometimes call the Jack function the {\it Jack polynomial} if $\bx$ is a finite sequence,
$\bx=(\bx_1,\dots,\bx_n)$ say. 
The Jack polynomials also satisfy another orthogonality condition.
Define the scalar product by
\begin{equation} 
\langle \phi, \psi \rangle_{n,\alpha}' = \frac{1}{n!} \int_{\bT^n}
\phi(z_1,\dots,z_n) \psi(z_1^{-1}, \dots, z_n^{-1}) |V(z_1,\dots,z_n)|^{2/\alpha} \,
\dd z_1 \cdots \dd z_n
\end{equation} 
for $n$-variables symmetric polynomials $\phi, \psi$.
Then we have the orthogonality
\begin{equation} \label{Eq:SecondOrthogonality}
\langle P_{\lambda}^{(\alpha)}, Q_{\mu}^{(\alpha)} \rangle_{n,\alpha}'
= \delta_{\lambda,\mu} I_n(\alpha) \prod_{(i,j) \in \lambda}
\frac{n+(j-1)\alpha-i+1}{n+j \alpha -i},
\end{equation}
where
\begin{equation} \label{Eq:ValueIn}
I_n(\alpha) := 
\frac{1}{n!} \int_{\bT^n}
|V(z_1,\dots,z_n)|^{2/\alpha} \,\dd z_1 \cdots \dd z_n
= \frac{\Gamma(n/\alpha+1)}{n! \, \Gamma(1/\alpha+1)^n}.
\end{equation}
The last equality is given in \cite[\S 8]{AAR} for example.

\subsection{Evaluation of Toeplitz hyperdeterminants in terms of Jack functions}

Let $\mathbf{1}$ be the function such that $\mathbf{1}(z)=1$ for any $z \in \bT$.
From Theorem \ref{Thm:HSF} and expression \eqref{Eq:ValueIn}, we have
$$
D_{n}^{[2m]}(\mathbf{1}) = I_n(1/m) = \frac{(mn)!}{n! \, (m!)^n}.
$$
For a function $f$ on $\bT$, denote by $\widehat{D}_n^{[2m]}(f)$ the normalized
Toeplitz hyperdeterminant 
$$
\widehat{D}_n^{[2m]}(f) = \frac{D_n^{[2m]}(f)}{D_{n}^{[2m]}(\mathbf{1})}.
$$

We compute Toeplitz hyperdeterminants employing Jack polynomials.
Let $f$ be a function in $L^1(\bT)$ with the Fourier expansion 
$f(z)= \sum_{k \in \bZ} d(k) z^k$.
The value $D^{[2m]}_n(f)$ is independent of $d(k)$ for $|k|>(n-1)m$
because the $d(k)$ do not appear among entries of $D^{[2m]}_n(f)$.
Thus, there exists a non-negative integer $R$ such that 
$D^{[2m]}_n(f)= D^{[2m]}_n(F_R)$,
where $F_{R}(z)= \sum_{ k \ge -R} d(k) z^k$.
Then, from Theorem \ref{Thm:HSF}, we see that
$$
D_{n}^{[2m]}(F_R) = \frac{1}{n!}
\int_{\bT^n} \prod_{k=1}^n z_k^R F_R(z_k) \cdot \overline{(z_1 \cdots z_n)^R} \cdot
|V(z_1,\dots,z_n)|^{2m} \, \dd z_1 \cdots \dd z_n.
$$
We have
$P_{(R^n)}^{(\alpha)}(\bx_1,\dots,\bx_n) = (\bx_1 \cdots \bx_n)^R$
for any $\alpha>0$ (see \cite[\S VI (4.17)]{Mac}), while
the formal power series 
$$
S_f(\bx_1,\dots,\bx_n;R) := \prod_{k=1}^n \bx_k^R F_R(\bx_k)
$$
in $\bC[[\bx_1,\dots, \bx_n]]$
can be expanded with respect to Jack $Q$-polynomials.
Therefore, if we obtain the coefficient of $Q_{(R^n)}^{(1/m)}(\bx_1,\dots,\bx_n)$
in this expansion, 
then we obtain the value $D_{n}^{[2m]}(F_R)$ 
by the orthogonality \eqref{Eq:SecondOrthogonality}.
Indeed,
if we denote by $\gamma(f,n,m,R)$ the coefficient of $Q_{(R^n)}^{(1/m)}$
in the expansion of $S_f(\bx_1,\dots,\bx_n;R)$ with respect to 
Jack $Q$-polynomials $Q_{\lambda}^{(1/m)}$,
we have
\begin{align*}
D_{n}^{[2m]}(f)
=& \frac{1}{n!} \int_{\bT^n} S_f(z_1,\dots,z_n;R) P_{(R^n)}^{(1/m)}(z_1^{-1},\dots,z_n^{-1})
|V(z_1,\dots,z_n)|^{2m} \dd z_1 \cdots \dd z_n \\=&  \gamma(f,n,m,R) \cdot 
\langle Q_{(R^n)}^{(1/m)}, P_{(R^n)}^{(1/m)} \rangle_{n,1/m}'.
\end{align*}
The explicit value of $\langle Q_{(R^n)}^{(1/m)}, P_{(R^n)}^{(1/m)} \rangle_{n,1/m}'$
can be obtained from \eqref{Eq:SecondOrthogonality}.
Finally, we have the following theorem.

\begin{thm} \label{Thm:Toeplitz}
Let $f \in L^1(\bT)$ and let $R$ be 
a non-negative integer such that $D^{[2m]}_n(f) = D^{[2m]}_n(F_R)$.
Then
we have
$$
\widehat{D}_{n}^{[2m]}(f) =\gamma(f,n,m,R)
\prod_{i=1}^n \prod_{j=1}^R \frac{i m +j-1}{(i-1)m+j},
$$
where $\gamma(f,n,m,R)$ is the coefficient of $Q_{(R^n)}^{(1/m)}$
in the expansion of $S_f(\bx_1,\dots,\bx_n;R)$ with respect to 
Jack $Q$-polynomials $Q_{\lambda}^{(1/m)}$. \qed
\end{thm}

As this theorem indicates,
the computation of a Toeplitz hyperdeterminant is reduced
to the evaluation of the coefficients of a Jack polynomial expansion.
However, it is hard to obtain explicit values of $\gamma(f,n,m,R)$ in general.
Here we give a few simple examples where $\gamma(f,n,m,1)$ can be explicitly calculated.

\begin{example}
Consider $f(z)=z^a-z^{-1}$
where $a$ is  a positive integer. 
In the notation of Theorem \ref{Thm:Toeplitz},
we can take $R=1$.
Then
$S_f(\bx_1,\dots,\bx_n;R) = (-1)^n \prod_{k=1}^n(1-\bx_{k}^{a+1})$.
The degree of each term in the polynomial $\prod_{k=1}^n (1-\bx_k^{a+1})$ 
is divisible by $a+1$  and therefore
$\gamma(f,n,m,R)=0$ unless ${\displaystyle n \equiv 0 \pmod{a+1}}$.
Henceforth, assume $n \equiv 0 \pmod{a+1}$ and put $n_a=n/(a+1)$. 
Then the term of degree $n$ in $(-1)^n \prod_{k=1}^n(1-\bx_{k}^{a+1})$
is given by an elementary symmetric polynomial
\begin{align*}
& (-1)^{n+n_a} e_{n_a}(\bx_1^{a+1},\dots,\bx_n^{a+1}) =
(-1)^{n+ n_a} (e_{n_a} \circ p_{a+1}) (\bx_1,\dots, \bx_n) \\
=& \sum_{\lambda \vdash n_a} \frac{(-1)^{n+\ell(\lambda)}}{ z_{\lambda}} 
(p_{\lambda} \circ p_{a+1})(\bx_1,\dots,\bx_n).
\end{align*}
Here $\circ$ denotes the plethysm product (\cite[\S I-8]{Mac})
and, in the second equality, the formula (2.14') in \cite[\S I]{Mac} is used.
Now we use the basic property $p_{\lambda} \circ p_k= p_{k\lambda}$, where 
$k\lambda=(k\lambda_1, k \lambda_2, \dots)$,
and the expansion formula (obtained from \cite[\S VI, (10.27)]{Mac})
\begin{equation} \label{Eq:Expand-p}
p_{\rho}= \alpha^{\ell(\rho)} z_{\rho} \sum_{\lambda : |\lambda|=|\rho|}
\frac{\theta_{\rho}^\lambda(\alpha)}{c_{\lambda}(\alpha)}Q_{\lambda}^{(\alpha)}.
\end{equation}
Then we have 
$\gamma(f,n,m,R) = 
\sum_{\lambda \vdash n_a} (-1)^{n+\ell(\lambda)} 
z_{\lambda}^{-1} m^{-\ell(\lambda)} z_{(a+1)\lambda}
\theta^{(1^n)}_{(a+1)\lambda}(\frac{1}{m}) / c_{(1^n)}(\frac{1}{m})$.
Finally, since $\theta_{\rho}^{(1^n)} (\alpha) = (-1)^{n-\ell(\rho)} n! \, z_{\rho}^{-1}$
(see \cite[\S VI-10, Ex.1]{Mac}) and 
$n! \sum_{\lambda \vdash n} \alpha^{\ell(\lambda)} z_{\lambda}^{-1} 
= \prod_{i=0}^{n-1}(\alpha +i)$,
we have 
$$
\gamma(f,n,m,R)= \frac{\prod_{i=0}^{n_a-1} (im+1)}
{(n_a)! \, m^{n_a}}.
$$
Therefore, by Theorem \ref{Thm:Toeplitz},
the Toeplitz hyperdeterminant 
of $f(z)=z^a-z^{-1}$ is given by
$$
\widehat{D}^{[2m]}_{n}(f)= \prod_{i=n_a}^{n-1} \frac{im+m}{im+1}
$$
if $n \equiv 0 \bmod a+1$, and $\widehat{D}^{[2m]}_n(f)=0$ otherwise.  \qed
\end{example}

\begin{example}
Let $f(z) =\{s z (1-sz) \}^{-1}$ with $0<|s|<1$.
Then ${\displaystyle S_f(\bx_1,\dots,\bx_n;R) = 
s^{-n} \prod_{k=1}^n (1-s \bx_k)^{-1}}$ with $R=1$.
The term of degree $n$ is the complete symmetric polynomial $h_n(\bx_1,\dots,\bx_n)$.
Hence, by a similar discussion to the last example,
we have 
$\gamma(f,n,m,R) = (-1)^{n-1} (n! \, m^n)^{-1} \prod_{i=1}^{n-1} (im-1)$.
Here the value is independent with $s$.
Therefore the Toeplitz hyperdeterminant is
$$
\widehat{D}^{[2m]}_{n}(f) = 
(-1)^{n-1} \prod_{i=1}^{n-1} \frac{i m-1}{im+1}.
$$
In particular, $D^{[2]}_{n}(f) =0$ for $n>1$. \qed
\end{example}

\begin{example}
Let $f(z)=z^{-1}e^z$. Then we have $S_f(\bx_1,\dots,\bx_n;1) = \exp(p_1(\bx))$
and therefore the term of degree $n$ is $p_{(1^n)}/n!$.
It follows by \eqref{Eq:Expand-p} that $\gamma(f,n,m,1) = (m^n n!)^{-1}$,
and so
the Toeplitz hyperdeterminant is
$$
\widehat{D}^{[2m]}_{n}(f) = \prod_{i=1}^{n-1} (im+1)^{-1}. \qed
$$
\end{example}

\subsection{Strong Szeg\"{o} limit theorem for Toeplitz hyperdeterminants}

We consider the asymptotic limit  of a Toeplitz hyperdeterminant of size $n$
in the limit as $n \to \infty$, whence
we obtain a hyperdeterminant analogue of the strong Szeg\"{o} limit theorem.

\begin{thm} \label{Thm:Szego}
Let $f(z) = \exp \(\sum_{k \in \bZ} c(k) z^k\)$ be a function on $\bT$
and assume 
\begin{equation} \label{AssumptionSzego}
\sum_{k \in \bZ} |c(k)| < \infty \quad \text{and} \quad
\sum_{k \in \bZ} |k| \, |c(k)|^2 < \infty.
\end{equation}
Then we have 
$$
\widehat{D}_n^{[2m]}(f) \sim \exp
\( c(0) n+ \frac{1}{m} \sum_{k=1}^\infty k c(k) c(-k) \) 
$$
as $n \to \infty$. \qed
\end{thm}

The case $m=1$ is actually the strong Szeg\"{o} limit theorem for a Toeplitz determinant,
see e.g. \cite{Mehta}.
This theorem follows from Theorem \ref{Thm:HSF} and the following lemma.
The lemma has previously been given in \cite{Johansson1988, Johansson1998},  while the present author has given
 a simpler algebraic proof in \cite{M}.

\begin{lem} 
Let $f$ be as in Theorem \ref{Thm:Szego}.
Then, for any $\alpha>0$, we have
$$
\lim_{n \to \infty} \frac{1}{I_n(\alpha) \, e^{n c(0)}}
\int_{\bT^n} \prod_{j=1}^n f(z_j) |V(z_1,\dots,z_n)|^{2/\alpha} \, \dd z_1 \cdots \dd z_n =
\exp \( \alpha \sum_{k=1}^\infty k c(k) c(-k) \). \qed
$$
\end{lem}

\begin{example}
Let $f(z) =e^{x(z-z^{-1})}$ with $x >0$.
The Fourier coefficients are Bessel functions, see e.g. \cite[Eq. (4.9.10)]{AAR}.
Then Theorem \ref{Thm:Szego} says
$\lim_{n \to \infty} \widehat{D}_{n}^{[2m]}(e^{x(z-z^{-1})}) = e^{-x^2/m}$. \qed
\end{example}

\begin{example}
Let $f(z)=(1+t z)^{w_1} (1+s z^{-1})^{w_2}$ with complex parameters $s,t,w_1,w_2$
satisfying $|s|, |t| < 1$.
Since $c(k)= w_1 (-1)^{k+1} t^k/k$ and $c(-k)= w_2 (-1)^{k+1} s^k/k$ for $k>0$, we have
$$
\lim_{n \to \infty} \widehat{D}^{[2m]}_{n}(f) 
= \exp \( \frac{w_1 w_2}{m} \sum_{k=1}^\infty \frac{(st)^k}{k} \)
=(1-st)^{-w_1 w_2/m}. \qed
$$
\end{example}

\begin{example}
Define $c(k)= |k|^{-1-x}$ if $k \not=0$ and $c(0)=0$ with $x >0$.
Then 
${\displaystyle \sum_{k \in \bZ} |c(k)| = 2 \zeta(1+x)}$ and 
$\sum_{k \in \bZ} |k| |c(k)|^2 = 2 \zeta(1+2x)$,
where $\zeta(s)$ is the Riemann zeta function $\zeta(s) = \sum_{n=1}^\infty n^{-s}$
for $\mathrm{Re} \, s >1$.
Therefore $f(z)= \exp \( \sum_{k=1}^\infty (z^k-z^{-k})/ k^{1+x} \)$
and we have 
${\displaystyle \lim_{n \to \infty} \widehat{D}^{[2m]}_{n}(f) = e^{\zeta(1+2x)/m}}$. \qed
\end{example}

%
\section{Jacobi-Trudi type formula for Jack functions of rectangular shapes}
%

The main result of the present study is presented in this section.
We obtain 
the Jacobi-Trudi formula for Jack functions of rectangular shapes by employing
 the Toeplitz hyperdeterminant studied in the previous section.

Put
$$
G_{\bx}^{(\alpha)}(z) = \prod_{i=1}^\infty (1-\bx_i z)^{-1/\alpha},
\qquad
E_\bx(z)= \prod_{i=1}^\infty (1+\bx_i z).
$$
These functions are the generating functions for one-row $Q$-functions and one-column $P$-functions,
$$
G_\bx^{(\alpha)}(z) = \sum_{r=0}^\infty Q_{(r)}^{(\alpha)}(\bx) z^r,
\qquad
E_\bx(z) = \sum_{r=0}^\infty P_{(1^r)}^{(\alpha)}(\bx) z^r.
$$
Put $g_r^{(\alpha)}=Q_{(r)}^{(\alpha)}$.
The function $P_{(1^r)}^{(\alpha)}$ is just equal to
the elementary symmetric function $e_{r}$.
We put $g^{(\alpha)}_r =e_r=0$ unless $r \ge 0$.

Consider a shifted Toeplitz hyperdeterminant defined by
$$
\widehat{D}_{n;a}^{[2m]}(f) = \widehat{D}_{n}^{[2m]}(z^{-a} f(z)), \qquad a \in \bZ. 
$$
Then we have the following formula.

\begin{thm} \label{Thm:JackToeplitzFormula}
Let $m$, $n$, and $L$ be positive integers,
and let $\bx$ be the sequence of infinitely many variables $\bx=(\bx_1,\bx_2, \dots)$.
Then we have
\begin{equation} \label{eq:JackToeplitzF}
Q_{(L^n)}^{(1/m)} (\bx)= \widehat{D}^{[2m]}_{n;L}(G_\bx^{(1/m)}), \qquad
P_{(n^L)}^{(m)}(\bx) = \widehat{D}^{[2m]}_{n;L} (E_\bx).
\end{equation}
In other words, we have
\begin{align*}
Q_{(L^n)}^{(1/m)}=& \frac{n! \, (m!)^n}{(mn)!} \cdot
\hdet{2m}(g^{(1/m)}_{L+i_1+ \cdots+i_m-i_{m+1}- \cdots -i_{2m}})_{[n]},
\\
P_{(n^L)}^{(m)} =& \frac{n! \, (m!)^n}{(mn)!} \cdot
\hdet{2m}(e_{L+i_1+ \cdots+i_m-i_{m+1}- \cdots -i_{2m}})_{[n]}.
\end{align*}
\end{thm}

\begin{proof}
We may assume each $\bx_j$ is a non-zero complex number and $|\bx_j|$ is sufficiently small.
We apply Theorem \ref{Thm:Toeplitz} to the function $f(z)=z^{-L} E_\bx(z)$ and $R=L$.  
Then in the notation of Theorem \ref{Thm:Toeplitz} 
we have $S_f(z_1,\dots,z_n;L)= \prod_{k=1}^n \prod_{i=1}^\infty (1+\bx_i z_k)$.
The dual Cauchy formula for Jack functions (\cite[\S VI-5,10]{Mac}) says that
$$
\sum_{\lambda: \ell(\lambda) \le n} Q_{\lambda'}^{(1/\alpha)}(\bx) 
Q_{\lambda}^{(\alpha)}(z_1,\dots,z_n) = \prod_{i = 1}^\infty \prod_{k=1}^n
(1+\bx_i z_k).
$$
Therefore we obtain 
$\gamma(f,n,m,L)= Q_{(n^L)}^{(m)}(\bx)$ and 
$$
\widehat{D}_{n;L}^{[2m]}(E_\bx)= \widehat{D}_{n}^{[2m]}(f)= 
Q_{(n^L)}^{(m)}(\bx) \prod_{i=1}^n \prod_{j=1}^L \frac{im+j-1}{(i-1)m+j}
=P_{(n^L)}^{(m)}(\bx),
$$ 
and so we have obtained the second in equation \eqref{eq:JackToeplitzF} to be proved.

Recall the endomorphism $\omega_{\alpha}$ on the $\bC$-algebra of symmetric functions
(see \cite[\S VI-10]{Mac}). 
It satisfies the duality $\omega_{\alpha}(P_{\lambda}^{(\alpha)})=Q_{\lambda'}^{(1/\alpha)},$
and  so the first formula follows from the second formula in equation \eqref{eq:JackToeplitzF}.
\end{proof}

As a corollary of the previous theorem, we see that
$$
Q^{(1/m)}_{(L^n)} \in \bQ[ g_{L+i}^{(1/m)} : -m(n-1) \le i \le m(n-1)], \qquad 
P^{(m)}_{(n^L)} \in \bQ[ e_{L+i} : -m(n-1) \le i \le m(n-1)].
$$

Since the Schur function is the Jack function associated with $\alpha=1$: $s_{\lambda}=
Q_{\lambda}^{(1)}=P_{\lambda}^{(1)}$,
the case $m=1$ in Theorem \ref{Thm:JackToeplitzFormula} reduces to the 
well-known Jacobi-Trudi identity
and the dual identity
for the Schur function of a rectangular shape:
$$
s_{(L^n)} = \det(h_{L-i+j})_{1 \le i,j \le n}, 
\qquad s_{(n^L)}  = \det(e_{L-i+j})_{1 \le i,j \le n}.
$$

If $\alpha$ is an even number and its inverse, 
the Jack function of a rectangular shape
can also be expressed by a hyperpfaffian as follows.
From Theorem  \ref{Thm:JackToeplitzFormula} and Theorem \ref{Thm:ToeplitzPfaffian},
we have the following corollary.

\begin{cor}
For any positive integers $m,L$ and $n$, we have 
\begin{align*}
Q^{(1/(2m))}_{(L^n)} =& \frac{n! ((2m)!)^n}{(2mn)!} 
\hpf{2m} \(\prod_{s=1}^m (i_{2s}-i_{2s-1}) \cdot  
g^{(1/(2m))}_{L+m(2n+1)-i_1- \cdots- i_{2m}}\)_{[2n]}, \\
P^{(2m)}_{(n^L)} =& \frac{n! ((2m)!)^n}{(2mn)!} 
\hpf{2m} \(\prod_{s=1}^m (i_{2s}-i_{2s-1}) \cdot  
e_{L+m(2n+1)-i_1- \cdots- i_{2m}}\)_{[2n]}.
\end{align*}
\qed
\end{cor}

\section{Appendix: Barvinok's hyperpfaffian} \label{Appendix}
%

Our hyperpfaffian $\hpf{2m}$ defined by expression \eqref{Eq:DefHyperPfaffian}
differs from
Barvinok's hyperpfaffian $\hPf{2m}$ in \cite{Barvinok}.
In this section, we give the explicit relationship between these two hyperpfaffians .

For positive integers $m$ and $n$,
we put
$$
\mathfrak{E}_{2m n,2m} = \{ \sigma \in \mf{S}_{2m n} \ | \ 
\sigma(2m (i-1) +1) < \sigma (2m (i-1) +2 ) < \cdots < \sigma ( 2m i) \ ( 1 \le i \le n )\}.
$$
Let $M=(M(i_1, \dots, i_{2m}))_{[2m n ]}$ 
be an array  satisfying 
$M(i_{\tau(1)}, \dots, i_{\tau(2m)}) = \sgn(\tau) M(i_1, \dots, i_{2m})$
for any $\tau \in \mf{S}_{2m}$.  
Barvinok \cite{Barvinok} defines his hyperpfaffian by 
$$
\hPf{2m}(M) = \frac{1}{n!} 
\sum_{\sigma \in \mf{E}_{2m n, 2m}} \sgn(\sigma) \prod_{i=1}^n 
M(\sigma (2m(i-1) +1), \sigma (2m(i-1) +2), \dots, \sigma(2m i)),
$$
see also \cite{LuqueThibon2002}.

As the following proposition states, 
our hyperpfaffian $\hpf{2m}$ is expressed by  $\hPf{2m}$.

\begin{prop}
Let $B=(B(i_1,\dots,i_{2m}))_{[2n]}$
be an array satisfying condition \eqref{EqAlter}.
Let 
$$
M=(M(i_1,\dots,i_{2m}))_{[2mn]}
$$ 
be the array
whose entries $M(i_1,\dots,i_{2m})$ are given as follows: if $i_1 < \cdots <i_{2m}$ and if
 there exist $1 \le r_1, \dots, r_{2m} \le 2n$
such that $i_{2s-1}=2n(s-1) +r_{2s-1}$ and $i_{2s}=2n(s-1)+r_{2s}$ for any $1 \le s \le m$,
then $M(i_1,\dots,i_{2m})=B(r_1,\dots,r_{2m})$. Otherwise
 define $M(i_1,\dots,i_{2m})=0$.  
Then we have $\hpf{2m}(B)=\hPf{2m}(M)$.
\end{prop}

\begin{proof}
The value
$\prod_{i=1}^n 
M(\sigma (2m(i-1) +1), \sigma (2m(i-1) +2), \dots, \sigma(2m i))$ is zero
unless the permutation $\sigma \in \mf{E}_{2mn,2m}$ satisfies
$2n(s-1)+1 \le \sigma(2m(i-1)+2s-1), \sigma(2m(i-1)+2s) \le 2ns$
for any $1 \le s \le m$ and $1 \le i \le n$.
Therefore
\begin{align*}
\hPf{2m}(M) =&
\frac{1}{n!} \sum_{\sigma_1,\dots, \sigma_{m} \in \mf{E}_{2n}}
\epsilon(\sigma_1,\dots, \sigma_m) \\
& \times \prod_{i=1}^n
M(\cdots \cdots, 2(s-1)n+\sigma_{s}(2i-1), 2(s-1)n+\sigma_{s} (2i), \cdots \cdots),
\end{align*}
where $\epsilon(\sigma_1,\dots, \sigma_m)$ is the signature of the permutation $\sigma$
defined by 
$$
\sigma(2m(i-1) + 2s-1) = 2n(s-1) +\sigma_s(2i-1) \quad \text{and} \quad
\sigma(2m(i-1)+2s) =  2n(s-1) +\sigma_s(2i)
$$
for any $1 \le i \le n$ and $1 \le s \le m$.
Hence we have
$\hPf{2m}(M)= \epsilon(\mathrm{id}, \dots, \mathrm{id}) \hpf{2m}(B)$,
and so $\epsilon(\mathrm{id}, \dots, \mathrm{id})= \sgn(\rho)$,
where 
$$
\rho(2m(i-1)+2s-1) = 2n(s-1) +2i-1 \quad \text{and} \quad 
\rho(2m(i-1)+2s) =  2n(s-1) +2i.
$$
Now it is straightforward to see $\sgn(\rho)=1$.
\end{proof}

\medskip
\medskip
\noindent
{\bf Acknowledgement.} 
I would like to thank Professor Masao Ishikawa and Professor Hiroyuki Tagawa
for their comments.
I also appreciate Reviewer's suggestions for the revision.


\end{document}